\documentclass[10pt,leqno]{amsart}

\usepackage{amssymb}
\usepackage{amsmath}
\usepackage{amsthm}
\usepackage{latexsym}
\usepackage{amscd, xypic}
\usepackage[dvipsnames]{xcolor}
\usepackage{tikz-cd} 
\tikzcdset{arrow style=tikz}
\usepackage{a4wide}
\usepackage{enumitem}
\setlist[itemize]{leftmargin=*}
\setlist[enumerate]{leftmargin=*}

\begin{document}

\newtheorem{thm}{Theorem}
\newtheorem{prop}[thm]{Proposition}
\newtheorem{lem}[thm]{Lemma}
\newtheorem{cor}[thm]{Corollary}
\newtheorem*{que}{Question}
\newtheorem*{ass}{Assumption}
\theoremstyle{definition}
\newtheorem{defi}[thm]{Definition}
\newtheorem{ex}[thm]{Example}
\newtheorem{rem}[thm]{Remark}

\newcommand{\ob}{\mbox{ob}}
\newcommand{\mor}{\mbox{mor}}
\newcommand{\iso}{\mbox{Iso}}
\newcommand{\id}{\mbox{id}}
\newcommand{\G}{\mathcal {G}}
\newcommand{\U}{\mathcal {U}}
\newcommand{\N}{\mathbb {N}}
\newcommand{\si}{\sigma}
\newcommand{\rh}{\rho}
\newcommand{\ta}{\tau}
\newcommand{\rmod}{$R$\mbox{-mod}}
\newcommand{\grmod}{\mbox{$G$-$R$-mod}}
\newcommand{\HOM}{{\rm HOM}}
\newcommand{\Hom}{\mbox{Hom}}
\newcommand{\END}{{\rm END}}
\newcommand{\Ker}{\mbox{Ker}}
\newcommand{\Image}{\mbox{Im}}
\newcommand{\m}{^{-1}}

\newcommand{\GR}{\mbox{$\G$-$R$}}
\newcommand{\AbG}{\mbox{$\mbox{Ab}_{\G}$}}
\newcommand{\Ab}{\mbox{Ab}}
\newcommand{\rmd}{\mbox{$R$-md}}
\newcommand{\mdr}{\mbox{md-$R$}}
\newcommand{\rmds}{\mbox{$R$-md-$S$}}
\newcommand{\smod}{\mbox{$S$-mod}}
\newcommand{\rmods}{\mbox{$R$-mod-$S$}}
\newcommand{\modr}{\mbox{mod-$R$}}
\newcommand{\gsmod}{\mbox{$\G$-$S$-mod}}
\newcommand{\gmodr}{\mbox{$\G$-mod-$R$}}
\newcommand{\grmods}{\mbox{$\G$-$R$-mod-$S$}}
\newcommand{\grmodr}{\mbox{$\G$-$R$-mod-$R$}}
\newcommand{\gmods}{\mbox{$\G$-mod-$S$}}
\newcommand{\grmd}{\mbox{$\G$-$R$-md}}
\newcommand{\gsmd}{\mbox{$\G$-$S$-md}}
\newcommand{\gmdr}{\mbox{$\G$-md-$R$}}
\newcommand{\gmds}{\mbox{$\G$-md-$S$}}
\newcommand{\grmds}{\mbox{$\G$-$R$-md-$S$}}

\newcommand{\rumod}{R\text{-\textup{\textbf{umod}}}}
\newcommand{\sumod}{S\text{-\textup{\textbf{umod}}}}
\newcommand{\modgr}{\mbox{{\bf mod}-}\G\mbox{-$R$}}
\newcommand{\grumod}{\G\textup{-}\rumod}
\newcommand{\gsumod}{\G\textup{-}\sumod}
\newcommand{\umodgr}{\G\text{-\textup{\textbf{umod}}}\textup{-}R}
\newcommand{\umodgs}{\G\text{-\textup{\textbf{umod}}}\textup{-}S}

\title{A short noncomplex proof for the solution set of homogeneous second order linear differential equations with constant coefficients}

\author[P. Lundstr\"{o}m]{Patrik Lundstr\"{o}m}

\address{University West,
Department of Engineering Science, 
SE-46186 Trollh\"{a}ttan, Sweden}
\email{{\scriptsize patrik.lundstrom@hv.se}
}



\begin{abstract}
We provide a short proof, not utilizing complex numbers, 
for the solution set of homogeneous second order 
linear differential equations with constant coefficients.
\end{abstract}

\maketitle

We argue that the solution set
of second order linear differential equations
with constant 
\linebreak
coefficients 
should be presented using hyperbolic
functions in the distinct real root case: 

\vspace{2mm}

\noindent
{\bf Theorem.}
Let $a,b \in \mathbb{R}$. Put
$\alpha = a/2$ and $\beta = \sqrt{|a^2/4-b|}$.
All solutions to 
$y'' + ay' + by = 0$ are given by
(1) $y = e^{-\alpha x} ( C_1 \cosh(\beta x) + 
C_2 \sinh(\beta x) )$, if $a^2/4 > b$,
(2) $y = e^{-\alpha x} ( C_1 + C_2 x )$, 
\linebreak
if $a^2/4 = b$, and
(3) $y = e^{-\alpha x}(C_1 \cos(\beta x) + 
C_2\sin(\beta x))$,
if $a^2/4 < b$,
for arbitrary $C_1,C_2 \in \mathbb{R}$.

\vspace{2mm}

There are many reasons for doing so.
First of all, this presentation
highlights the beautiful symmetry between (1) and (3).
Secondly, if we in (1) and (3) use the 
first order approximations $\cosh(\beta x) \approx 1$,
$\sinh(\beta x) \approx \beta x$,
$\cos(\beta x) \approx 1$ and 
$\sin(\beta x) \approx \beta x$, for small $\beta$, 
then (2) can be seen, qualitatively, 
as a limiting case of both (1) and (3), as $\beta \to 0$.
Thirdly, and perhaps most importantly,
utilizing a condensed version of a clever trick by 
Dobos \cite{dobos2007}, and here expanding
it to to the real case also, 
one can construct
a short coherent elementary 
noncomplex proof of (1)-(3).

Namely, it is a straighforward calculation to show that the expressions 
in (1), (2) and (3) indeed are 
solutions to the differential equation.
Now suppose that $y$ is an arbitrary solution and put 
\linebreak
$z = e^{\alpha x} y$.
Then $z' = \alpha z + e^{\alpha x} y'$ and thus
$z'' = \alpha^2 z + 2 \alpha e^{\alpha x} y' + e^{\alpha x}y''
= e^{\alpha x} ((a^2/4)y + a y' + y'')$. 

(1) Let $a^2/4 > b$. Then  
$z'' - \beta^2 z = 
e^{\alpha x} ( (a^2/4)y + a y' + y'' - (a^2/4-b)y ) =
e^{\alpha x} ( y'' + ay' + by ) = 0$ 
and hence
$
( - z \sinh(\beta x) + z' \cosh(\beta x)/\beta )' = 0$
and
$
(-z' \sinh(\beta x)/\beta + z \cosh(\beta x))' =
0
$
which implies that
$
-z \sinh(\beta x) + z' \cosh(\beta x)/\beta = C_2
$
and
$
-z' \sinh(\beta x)/\beta + z \cosh(\beta x) = C_1
$
for some $C_1 , C_2 \in \mathbb{R}$.
Therefore
$
C_1 \cosh(\beta x) + C_2 \sinh(\beta x) = 
z \cosh^2(\beta x) - z \sinh^2(\beta x) = z
= e^{\alpha x} y
$
and thus finally we get that
$y = e^{-\alpha x} z = e^{-\alpha x} ( C_1 \cosh(\beta x) + 
C_2 \sinh(\beta x) )$.

(2) Suppose that $a^2/4 = b$. Then
$z'' = e^{\alpha x} ((a^2/4)y + a y' + y'') = 
e^{\alpha x} ( y'' + ay' + by ) = 0$. 
\linebreak
Hence
$z = C_1 + C_2 x$ for some $C_1,C_2 \in \mathbb{R}$
and therefore we get that
$y = e^{-\alpha x} z = e^{-\alpha x}(C_1 + C_2 x)$.

(3) Let $a^2/4 < b$. Then 
$z'' + \beta^2 z = e^{\alpha x} 
( (a^2/4)y + a y' + y'' + (b - a^2/4)y ) =
e^{\alpha x} ( y'' + ay' + by ) = 0$
and hence  
$
( z \sin(\beta x) + z' \cos(\beta x)/\beta )' = 0
$
and
$(-z' \sin(\beta x)/\beta + z \cos(\beta x))' = 0$
which implies that
\linebreak
$
z \sin(\beta x) + z' \cos(\beta x)/\beta = C_2
$
and
$-z' \sin(\beta x)/\beta + z \cos(\beta x) = C_1$
for some $C_1,C_2 \in \mathbb{R}$.
\linebreak
Therefore
$
C_1 \cos(\beta x) + C_2 \sin(\beta x) = 
z \cos^2(\beta x) + z \sin^2(\beta x) = z
= e^{\alpha x} y
$
and thus finally we get that
$y = e^{-\alpha x} z = e^{-\alpha x} ( C_1 \cos(\beta x) + 
C_2 \sin(\beta x) )$.
\qed

\vspace{2mm}

Note that the solutions in (1) 
above 
of course easily can 
be rewritten in the familiar form
\linebreak
$y = e^{-\alpha x} ( C_1 \cosh(\beta x) + 
C_2 \sinh(\beta x) ) = 
D_1 e^{r_1 x} + D_2 e^{r_2 x} $ where 
$D_1 = (C_1+C_2)/2$, $D_2 = (C_1-C_2)/2$ and
$r_1 = -\alpha + \beta$ and $r_2 = -\alpha - \beta$ 
are the two roots of the characteristic equation
$r^2 + ar + b = 0$.

\end{document}